\input amstex
 \documentstyle{amsppt}
 \magnification=\magstep1
 \vsize=24.2true cm
 \hsize=15.3true cm
 \nopagenumbers\topskip=1truecm
 \headline={\tenrm\hfil\folio\hfil}

 \TagsOnRight

\hyphenation{auto-mor-phism auto-mor-phisms co-homo-log-i-cal co-homo-logy
co-homo-logous dual-izing pre-dual-izing geo-metric geo-metries geo-metry
half-space homeo-mor-phic homeo-mor-phism homeo-mor-phisms homo-log-i-cal
homo-logy homo-logous homo-mor-phism homo-mor-phisms hyper-plane hyper-planes
hyper-sur-face hyper-sur-faces idem-potent iso-mor-phism iso-mor-phisms
multi-plic-a-tion nil-potent poly-nomial priori rami-fication sin-gu-lar-ities
sub-vari-eties sub-vari-ety trans-form-a-tion trans-form-a-tions Castel-nuovo
Enri-ques Lo-ba-chev-sky Theo-rem Za-ni-chelli in-vo-lu-tion Na-ra-sim-han Bohr-Som-mer-feld}

\define\rest#1{_{\textstyle{|}#1}} 

\define\Span#1{\left<#1\right>} 

\define\half{{\textstyle{1\over2}}}


\define\C{\Bbb C} 
\define\R{\Bbb R} 
\define\Z{\Bbb Z} 

\define\proj{\Bbb P} 

\define\sA{{\Cal A}} 
\define\sM{{\Cal M}} 

\define\al{\alpha}

\define\de{\delta}

\define\ga{\gamma}

\define\om{\omega}
\define\si{\sigma}
\define\De{\Delta}

\define\La{\Lambda}
\define\Om{\Omega}

\define\Si{\Sigma}








 \document

 \topmatter
  \title Fano versus Calabi - Yau \endtitle
  \author    Andrey Tyurin\footnote{The last lecture
of Professor A.N. Tyurin (24.02.1940 - 27.10.2002)
given at the Fano conference (Turin, October 2002).
The text was prepared and edited by Nik. Tyurin 
and Yu. Tiourina. We'd like to thank Max - Planck
- Institute for Mathematics (Bonn, Germany).} 
                \endauthor

   \address Steklov Institute (Moscow)
   \endaddress
  \email ntyurin\@thsun1.jinr.ru
      jtyurin\@mpim-bonn.mpg.de
   \endemail

\abstract  In this article we discuss some numerical parts
of the mirror conjecture. For any 3 - dimensional
Calabi - Yau manifold author introduces a generalization
of the Casson invariant known in 3 - dimensional geometry,
which is called Casson - Donaldson invariant. In the 
framework of the mirror relationship it  corresponds
to the number of SpLag cycles which are Bohr - Sommerfeld
with respect to the given polarization. To compute
the Casson - Donaldson invariant the author uses well known
in classical algebraic geometry degeneration principle. By it, when
the given Calabi - Yau manifold is deformed to a pair
of quasi Fano manifolds glued upon some K3 - surface,
one can compute the invariant in terms of "flag
geometry" of the pairs (quasi Fano, K3 - surface). 
\endabstract

   \endtopmatter

\head Introduction  \endhead

According to the "Oxford program" proposed by S. Donaldson and R. Thomas
([4]), some constructions of real gauge theories can be repeated after
generalization to the complex case. In a sense, it is a variant of
Arnold sentence which claims that every real notion has (or should have)
an analogy in  the complex case ([1]). So it is a kind
of complexification. The authors of [4] propose two possible ways allowing
to translate the Chern - Simons theory and the Yang - Mills theory
originally defined over 3 - dimensional and 4 - dimensional
real manifolds respectively to the cases of 3 - dimensional and 4
- dimensional complex manifolds. The first move comes with 
the question what should substitute the notion of orientation in the complex case.
Reasoning in the way parallel to the original real case, authors of [4]
show that compact complex manifold is c - orientable iff
its canonical class is trivial. Thus they restrict the investigation
by the condition which imposes that we deal with in the case of Calabi - Yau
manifolds. Thus a c - orientation is defined by an appropriate trivialization of
the canonical class. Here we study the 3 - dimensional complex case which
contains 3 - dimensional Calabi - Yau manifolds. 

This case is of an enormous interest because of some physical applications:
the Calabi - Yau realm produces  collections of non negative numbers (so called
Casson - Donaldson invariant) which enter in the computations of physical theories.
The combinatorical structure of these collections is the following: the set
$\Cal C \Cal Y_3$ of all possible topological types of compact non singular
Calabi - Yau threefolds can be coded by the set of lattices (so called
Mukai lattices) equipped with bilinear forms (so called Mukai forms);
any pair $(M, \chi)$ of the objects represents the corresponding topological
K - algebra of the vector bundles (in the broad sense of the reflexive sheaves)
over a given Calabi - Yau threefold; and the Casson - Donaldson invariant
is a map
$$
CD: \{(M, \chi)\} \to \Z,
$$
which gives us a collection of numbers derived from the given topological type.
Thus one of the problems posed by physics  is to compute 
that numbers which is a very hard task in general.

The original real construction gives a hint about  what can be exploited to perform
the computations in some special cases when the given Calabi - Yau manifolds
are constructive. It means that for a given $M$ there exists its birational
model such that a deformation to a reducible threefold is possible:
$$
M_0 = F_+ \cup_S F_-,
\tag 0.1
$$
where $F_{\pm}$ are two quasi Fano components, glued together transversally
along a K3 - surface $S$ which belongs simultaneously to both the anticanonical
systems of the components. In this setup:

1) the Mukai lattice of M can be reproduced from the Mukai lattices of
the components $F_{\pm}$

and

2) CD - invariant can be derived from the geometry of the Fano
varieties.

Thus in the constructive case the knowledge of  the Fano geometry
is sufficient to perform the desired computations.

Here we have to note that as it was emphasized by the authors of [4]
the collection of CD - invariants is an analogy of the classical Casson invariant
of compact oriented real 3 - manifolds. Since the classical one is computed
using the Heegard decomposition, in the complex case one could exploit the same idea.
As an intermediate step in the decomposition for the real manifolds one considers 
the following picture: a Riemannian surface $\Si$ is included into a given
3 - dimensional manifold $X$, cutting it into two pieces glued over the surface. Thus the picture is
$$
X = H_+ \cup_{\Si} H_-,
$$
where $H_{\pm}$ is an oriented  3 - manifold with boundary such that
$$
\partial H_{\pm} = \pm \Si.
$$
This is exactly the same as in (0.1) since the complex analogy of real 3 - dimensional
oriented manifold with boundary is quasi Fano variety with a fixed K3 - surface
from the anticanonical system. This K3 - surface plays the role of the complex
boundary: since it should be c - oriented it is a K3 - surface. 

This parallel imposes the following sentence: since every real 3 - dimensional
oriented compact manifold can be represented
by an appropriate Heegard diagram, so we can expect (basing on our practical
experience) that every Calabi - Yau manifold in dimension 3 is
constructive and thus we can restrict the investigations of non linear
super $\si$ - models with CY - target spaces to the Fano realm. Also,
if it is true it should be just a finite set of different topological
types for 3 - dimensional Calabi - Yau manifolds.

\head 1. Mirror picture: LHS and RHS
\endhead

As it was observed first by theoretical physicists,
 for some Calabi - Yau 3 - manifolds one has so - called
mirror partners  such that if $X$ and $X'$
are in the "mirror" correspondence then the deformation properties of
the first one correspond to some properties of the cubic intersection form
restricted to the Picard lattice over the second one ([14]). It implies in particular that the Hodge
diamonds of $X$ and $X'$ are  related by some reflection.  
Even more rough description recognizes $X$ and $X'$ as  mirror partners
if the equality holds only for two entries of the Hodge diamonds
namely
$$
h^{1,1}(X) = h^{1,2}(X')
$$
and vice versa. While this weakening of the mirror conditions leads to
pure numerical coincidences only, we shall understand the mirror correspondence
in the following way: two 3 - dimensional Calabi - Yau manifolds are
related by the mirror duality if the algebraic geometry over the first one
is isomorphic to the symplectic geometry over the second one and vice versa.
Of course it isn't quite clear how to treat both the geometries 
so one should explain first what do they mean as well as
define what kind of equivalence  is desired to claim that these geometries
are isomorphic. A slightly abstract approach to this problem is given by
so called homological mirror symmetry proposed by M. Kontzevich which
suggests some equivalence between two derivations from the geometries:
on the LHS (algebraic geometry) one takes the derived category of
coherent sheaves while on the RHS (symplectic geometry) it should be
the Fukaya category --- and these two are isomorphic as categories.
On the other hand, the algebraic geometry contains some "real"
objects, which have been constructed in boundaries of  its framework and
which are well understood and more than habitual to
a number of the world experts who gathered at our conference.
We mean the moduli spaces of vector bundles and sheaves and
the systems of submanifolds (f.e. complete linear systems
or related objects such as the Chow groups etc.). At the same time
one could expect that the RHS of the mirror picture admits
some objects of the same reality --- a moduli space 
and systems of some special submanifolds. 
From this point of view the RHS used to be "terra incognita"
up to the middle of nineties. At that time, different ways of
developing the "classical" symplectic geometry 
were consolidated and united such that the Floer
homology of Lagragian submanifolds, pseudoholomorphic curves (in dimension 4),
special Lagrangian cycles (for the Calabi - Yau case) etc. turned to be the parts of some unified
consistent theory.  Thus now one has the ingredients at the RHS which  can be
compared to the known objects  existing on the LHS. In this section
we briefly recall  the generating objects for some
numerical derivations on both the sides of the mirror picture.       

\subhead 2.1. LHS
\endsubhead

Let $M$ be a smooth complete algebraic Calabi - Yau threefold which is equipped with
a Ricci flat Kahler metric $g$, giving the corresponding Kahler form $\om$.
At the begining point, one fixes an additional data: a complex orientation
---  a choice of a holomorphic  (3,0) - form $\Om$
which is defined up to the phase scaling $e^{i \phi}$ since the norm is fixed
by the metric. We collect these data in the quadruple
$$
(M, g, \om, \Om).
$$
The ring of even dimensional cohomology of M 
$$
H^{2*}(M, \Z) = \oplus_{i=0}^3 H^{2i}(M, \Z)
$$
can be equipped with an involution $*$ acting componentwise by the following
formulas:
$$
\aligned
*|_{H^0(M, \Z) \oplus H^4(M, \Z)} = id, \\
*|_{H^2(M, \Z) \oplus H^6(M, \Z)} = - id.\\
\endaligned
$$
The definition implies that the bilinear form
$$
(u, v) = - [v^* \cdot u]_6
$$
is skew - symmetric (one could recognize the reason to introduce $*$ exactly in
the derivation of this skew - symmetric property). At the same time
one has a similar natural involution on the algebraic K - functor
$K^0_{alg}$ on $M$ which sends each vector bundle $E$ to the dual bundle $E^*$.
It's easy to see that these two involutions are related by the homomorphism
$$
ch: K^0_{alg} \to H^{2*}(M, \Bbb Q)
$$
such that the diagram
$$
\matrix  K^0_{alg} & @>ch>> & H^{2*}(M, \Bbb Q)\\
\downarrow^* &  & \downarrow^* \\ K^0_{alg} & @>ch>> & H^{2*}(M, \Bbb Q) \\
\endmatrix
$$
commutes (here $ch$ is the Chern character map).

On the other hand one has on $K^0_{alg}$ the following bilinear form
$$ 
- \chi(E_1, E_2) = \sum_{i=0}^3 (-1)^{i+1} rk Ext^i(M, E_1^* \otimes E_2),
$$
where the homological spaces are the coherent cohomology of sheaves. This form 
can be represented as the image of a bilinear form over $H^{2*}(M, \Z)$;
recall that by the Riemann - Roch - Hirzebruch theorem we have
$$
\chi(E_1, E_2) = [ch E^2 \cdot ch E^*_1 \cdot td_M]_6,
$$
where $td_M$ is the Todd class of $M$.  In
the special situation of the Calabi - Yau manifolds   bilinear form $\chi$
is skew symmetric:
$$
\aligned
H^0(M, E_1^* \otimes E_2)^* = H^3(E_2^* \otimes E_1) \\
H^1 (E_1^* \otimes E_2)^* = H^2(E_2^* \otimes E_1)\\
H^2(E_1^* \otimes E_2)^* = H^1(E_2^* \otimes E_1)\\
H^3(E^*_1 \otimes E_2)^* = H^0(E_2^* \otimes E_1)\\
 \endaligned
$$
by the Serre duality which in this case reads as
$$
H^i(E^* \otimes F) = H^{3-i}(F^* \otimes E)
$$
since the canonical class is trivial.
Moreover, the Todd class of $M$ equals
$$
td_M = 1 + \frac{1}{12}c_2(M)
$$
and consequently
$$
td_M^* = td_M
$$
since $c_2(M) = c_2^*(M)$ by the definition. Thus we have a special class
in $H^{2*}(M, \Bbb Q)$ namely
$$
\sqrt{td_M} = 1 + \frac{1}{24} c_2(M)
$$
defined uniquely by the condition that the leading term is 1. 

Following S. Mukai we slightly correct the Chern character map
twisting it by the element
$$
m(E) = ch E \cdot \sqrt{td_M} \in H^{2*}(M, \Z)
$$
for a vector bundle $E$; we call it the Mukai vector of $E$.
The  Mukai vectors form a lattice $L_M$ which is called
the Mukai lattice, equipped with the bilinear skew symmetric form
$\chi$; thus the topological type of the underlayng Calabi - Yau
3 - manifold can be encoded by the pair $(L_M, \chi)$.
 
On the other hand, if $E$ is a holomorphic bundle over $M$, then
one has the following derivations from the deformation theory:

1) the space $Ext^1 (E, E)$ is the space
of infinitesimal deformations of $E$;

2) the space $Ext^2 (E, E)$ is the space of obstructions;

and 

3) locally the moduli space of holomorphic vector bundles of
type $m(E)$ is presented by the Kuranishi map
$$
k: Ext^1 (E, E) \to Ext^2 (E, E)
$$
such that
$$
\Cal M_E|_{O(E)} = k^{-1}(0)
$$
in a small neighborhood $O(E)$ of the given point $E \in \Cal M_E$.

 But in the 
Calabi - Yau realm 
$$
Ext^1(E, E) = Ext^2(E, E)^*
$$
by the Serre duality. Thus one could expect in general that
the moduli spaces are zero dimensional. We claim here that
at least the virtual (expected) dimension of the moduli space of simple
vector bundles of a topological type $m(E) \in L_M$
is zero.

Further, the Kahler form $\om$ fixed on $M$ determines a polarization
$H$ on $M$ which separates special class of holomorphic vector bundles
--- the bundles which are stable with respect to the given polarization.
The stability condition imposes one additional requirement on the topological
types: for a vector bundle (or reflexive sheaf) $E$ the expression
$$
c_2(E) - \frac{rk(E) - 1}{2 \cdot rk(E)} \cdot c_1^2(E) = \De (E)
$$
is called the discriminant of the bundle. Then the Bogomolov inequality
necessary for any stable bundle is given by the constraint
$$
\De(E) \cdot H > 0
$$
for the polarization $H$ over $M$.

Roughly the idea of the Casson - Donaldson invariant is as follows:
the moduli space of $H$ - stable holomorphic vector bundles
of a given type $m \in H^{2*}(M, \Z)$ is expected to be a finite set
$$
\Cal M_H(m) = \{p_1, ..., p_k\},
$$
Thus, counting the number, one can define an integer function on the Mukai lattice
$$
CD_H(m) = deg \Cal M_H(m) \in \Z
$$
and call it the Casson - Donaldson invariant.
Of course it is an ideal picture: for example, we added
the stability condition since otherwise the set of holomorphic
vector bundles of a fixed type doesn't in general posess
any good structure. For stable bundles of a fixed topological type
$m \in L_M$ the moduli space admits  a natural compactification
$\overline{\Cal M^s_H(m)}$ which in the transversal case 
carries the structure of zero dimensional scheme. As the length of the scheme
the number
$$
CD_H(m) = deg \Cal M^s_H(m)
\tag 1.1
$$
is well defined. It is quite reasonable to call this number the Casson - Donaldson
invariant: on the one hand, it is obviously analogous to the Casson invariant
as it was proposed in [4]; on the other hand, it is an obvious analogy of the Donaldson polynomial of degree zero for real 4 - manifolds ([3]).

The problem is that in the real life the transversal situation
can not occur for some bundles: deformations can be unobstructed
what happens, for example,  in the case of the tangent bundle of $M$.
Then the geometrical dimension of the moduli space of the topological type 
of $TM$ is equal to the dimension of the component of $CY_3$ - moduli space. To 
get the number in this case we have to apply the usual "deformation
to the normal cone" trick.

The difference between the virtual and the infinitesimal dimensions can be measured as follows. Let us consider on the set of all vector bundles some other bilinear form
$h$ defined as
$$
h(E_1, E_2) = rk H^1(E_1^* \otimes E_2) - rk H^0(E_1^* \otimes E_2).
\tag 1.2
$$
Then it is easy to see that $\chi$ is exactly the skew symmetric part of
$h$. At the same time the symmetrical part of $h$ could be viewed as an analogy
of the sum of the Betti numbers. We use it instead of the usual Euler 
characteristics of coherent sheaves taking in mind the real case where
the sum of the Betti numbers gives a numerical estimation for the statement
of the Arnold's conjecture, instead of pure topological Euler characteristics
of the based real manifold. As in the symplectic geometry, in our case
the number given by $h$ is not a topological invariant. Now we see that 
the difference between the virtual and the infinitesimal
dimensions is measured by the bilinear form $h$. It is fruitful to consider
the vector bundles for which these dimensions coincide. Such simple bundles
are called exeptional, aspheric or spheric depending on the context (we will
see at the next subsection that they are analogous to some spheres
on the RHS). We will call such bundles exeptional because they are simple and
infinitesimally rigid. S. Mukai noted that the bundles play the role of
roots in the Mukai lattice. In particular, for every such a bundle
one has
$$
h(E, E) = -1
$$
which is the minimal number for simple bundles. At the same time the Casson
- Donaldson invariant is well defined so these objects are
"right" from the point of view of theoretical physics. One can really
derive a collection of numbers for Mukai vectors corresponding
to these exeptional bundles such that $CD_H(m)$ is the degree
(1.1) above.

As we've mentioned, the bilinear form (1.2)
is not purely topological. One could answer 
the following natural question about it:
is there an equivalence relation $\sim$ on $K^0_{alg}$ such 
that the form $h$ is equal to the lifting of a bilinear form
on the lattice $K^0_{alg}/ \sim = F$. Obviously, $F$ can not
be equal to $H^{2*}(M, \Z)$ and if the answer to the question posed
here is "Yes" then it should be a lattice $F$ between the Chow
ring $CH^* (M) = A^*(M)/(rational \quad equivalence)$ and $H^{2*}(M, \Z)$
(for example, $\Cal A(M) = A^*(M)/(algebraic \quad equivalence)$). At the same time
by the Griffits theorem for 3 - dimensional quintic $M \subset \C \proj^4$
we have $\Cal A(M) \neq H^{2*}(M, \Z)$, thus our question is of the "Lefshetz
problem" type.

 Further, recall that coherent sheaves and bundles over any variety can be transformed 
by so called modular operations: universal extensions,
universal divisions and some of their combinations ([11], Ch. 2).
These operations preserve the properties
of the moduli spaces. The regular application of these transformations
have begun with S. Mukai paper [8] where the author described the "reflection" operation which can be decomposed 
into a combination
of the divisions and the extensions (or "returns" 
in the special terminology used by the experts). The spectrum of the possible applications
of these operations is quite wide: from the derived categories of coherent sheaves to the categories of the representations 
of algebras and quivers.
We can extract some combinatorical part from the purely
geometric investigations using such operation (and almost always
we obtain some representation of a braid group). And this usually gives
some mutual underlayng rules for a priori different theories.

 Every Mukai vector $m \in L_M$ defines a transformation of the full lattice
given by the natural formula
$$
\al_m(m') = - m' - \chi(m, m') \cdot m,
$$
but we prefer to define the following corrected version
$$
\al_m (m') = - m' - h(m, m') \cdot m.
$$
From the first viewpoint it seems to be ill defined, because the form $h$ (1.2)
doesn't descend to the even cohomology ring. However, extending
the discussion to the derived category $D^b(M)$ of coherent sheaves
on $M$ one could go this way: let us consider the map to the Atiyah ring
$$
r: D^b(M) \to K^0_{alg}.
$$
then if $m$ can be realised by an exeptional bundle $E$ it can be shown that the desired transformation can be defined as a functor 
$$
\al_E: D^b(M) \to D^b(M)
$$
on the derived category. Of course, such  lifting heavily depends
on the choice of this exeptional bundle $E$
in the fixed topological class $m(E)$. Let us choose and fix it. Then for any 
bundle $E'$ we have the following sheaf
$$
\al_E(E') = ker (can: H^0(E^* \otimes E') \otimes E \to E').
$$
If this one is a stable bundle:
$$
[\al_E(E')] = \al_{[E]}([E']),
$$
then this bundle is the result of a modular operation, namely of the universal
division. 

Let $H \in Pic M$ be a polarization of $M$. This one as a bundle of rank
1 defines an automorphism of the Mukai lattice:
$$
[T_{H^k}(E)] = T^k_H = [E \otimes H^k],
$$
which can be lifted to the derived category or to any other "category"
which could be defined by vector bundles. Using the standard homological 
technique and wordwise the same arguments as in Ch. 2 of [11], one can easily
prove the following
\proclaim{Proposition (1.1)} For any two topological types $m$ and $m'$
with holomorphic realizations $\{E_1, ..., E_{CD_H(m)}\}$ and
$\{E_1', ..., E_{CD_H(m')}'\}$ by exeptional stable bundles there exists a level $k_0(m, m') \in \Bbb N$ such that

1) for any $k > k_0(m, m')$ all modular operations $\al_{E_i}(E_j' \otimes H^k)$ (??) are correctly defined;

2) all the bundles from the collection
$\{\al_{E_i}(E_j') \otimes H^k)\}$ are exeptional and stable   and
form the complete collection of the holomorphic realizations of a
topological class $\al_m(T^k_H(m'))$;

3) consequently, the Casson - Donaldson number for the last topological type is given by the formula 
 $$
CD_H(\al_m(T^k_H(m'))) = CD_H(m) \cdot CD_H(m').
$$
\endproclaim

From this observation one immediately gets
\proclaim{Corollary 1.1} The Casson - Donaldson invariant is
unbounded as a function.
\endproclaim

Of course, there are special vector bundles for which the number $CD_H$ is equal to 1. For example, 
$$
CD_H(m(L)) = 1
$$
for any line bundle $L$ over $M$ (recall that by the definition $h^{1,0} = 0$).
On the other hand, there exist many types of sheaves with special Mukai vectors whose moduli spaces 
are compact and smooth (but of positive dimensions). For example, take
$F = \Cal O_p$ where $p \in M$ is a point. In this case the deformation
to the normal cone gives us
$$
CD_H(m(\Cal O_p)) = \chi_{top}(M) = \sum_i (-1)^i rk H^i(M, \R)
$$
where the right hand side is just the topological Euler characteristics. As well
the interpretation  of the moduli space $\Cal M(m)$ as the zero set of a holomorphic
differential on the space $\Cal D_M''(E)$ of $\bar \partial$ - connections
(see [13]) shows that
$$
CD (m(E)) = \chi_{top}(\Cal D_M''(E))
$$
(in spite of the fact that the space $\Cal D_M'(E)$ is infinite dimensional itself). 
This space depends on the topological type of $E$ only, that is, on the Chern
character $ch(E)$ or on the Mukai vector $m(E)$.

At the end of the LHS short description we have to mention that obviously the Casson - Donaldson 
invariant is the imitation of a well known one: let $m \in H^2(M, \Z)$ and
$R(m)$ be the family of the rational curves which represent this cohomology class.
Then by the same reason from the deformation theory we can expect that scheme
$R(m)$ is zero dimensional of length
$$
R_M(m) = deg R(m).
$$
The function 
$$
R_M: H^2(M, \Z) \to \Z
$$
underlies the physical parameters in the same way as $CD_H$ from (1.1). For example,
if $M \subset \C \proj^4$ is a quintic and $m \in H^2(M, \Z)$ is the class of the projective line in $\C \proj^4$, 
then $R_M(m) = 2875$. Moreover, the Clemens conjecture
predicts that our expectations of "finitness" are true for a generic quintic.

We propose here the following conjecture closely analogous to the Clemens' one, namely:
\proclaim{Conjecture} Over generic quintic in $\C \proj^4$ every stable rank 2
vector bundle $E$ is infinitesimally rigid, that is, $H^1(ad E) = 0$.
\endproclaim

\subhead 1.2. RHS
\endsubhead

On the symplectic side of the mirror picture one deals with the same type 
  quadruple $(M, g, \om, \Om)$ understanding it first as a symplectic
manifold with integer symplectic structure $\om$. The integrability condition
means that there is a line bundle $H$ with the first Chern class
$$
c_1(H) = [\om]
\tag 1.3
$$
equipped with a hermitian connection $a \in \sA_h(H)$ such that
$$
F_a = 2 \pi i \om
\tag 1.4
$$
and the last condition defines this connection almost uniquely
up to the gauge transformations. 

In the symplectic setup we consider Lagrangian cycles as the geometric objects 
instead of vector bundles (the term "cycle" used instead of "submanifold" is discussed in [6]). 
Like for the vector bundles, one can impose some conditions on 
the Lagrangian cycles to get some appropriate objects (moduli spaces).
In our setup (the case of so - called fixed complex polarization) 
there are two basic definitions which distinguish from the space $\Cal L$ 
of all Lagrangian cycles two classes. We recall briefly what they are.
The first condition is more universal: it can be imposed in general situation when
the symplectic form is integer. Namely one calls a Lagrangian submanifold $S$
Bohr - Sommerfeld if the restriction to it of the prequantization pair 
$(H, a)$ defined by (1.3) and (1.4) admits a covariantly constant section
([6]). Indeed, since $S$ is Lagrangian and the curvature form of $a$ is proportional
to $\om$, the restriction to any Lagrangian submanifold should be
isomorphic to a pair (trivial bundle, flat connection). But if this flat connection
is gauge equivalent to the trivial product connection then one says that
$S$ is Bohr - Sommerfeld. In [6] one constructs some moduli spaces of
such Bohr - Sommerfeld Lagrangian cycles.

Besides fixing some hermitian structure on $H$ to define the Bohr - Sommerfeld
condition on the Lagrangian cycles, the definition of SpLag cycles needs some
complex polarization as the ruling ingredient. Generally  complex polarization
means that some appropriate integrable complex structure compatible
with the given symplectic structure is fixed over the base manifold. At the same time
over our Calabi - Yau threefold the complex structure is defined by
the choice of $\Om$ which was made at the begining of our discussion.
Thus for the oriented Calabi - Yau threefold $(M, \Om)$ special Lagrangian cycle
(spLag cycle for short) is a 3 - dimensional Lagrangian submanifold $S$ such that
the restriction $\Om|_S$ satifies
$$
Im \Om|_S \cong 0 \quad \quad {\text or} \quad \quad Re \Om_S \cong Vol_g(S),
$$
and these two conditions are equivalent ([12]). On the other hand,
there is another way to define the cycles based on considerations of
Gauss map and Gauss vector field over the space of all Lagrangian cycles
([12]); in the same paper one finds some intermediate condition
on Lagrangian cycles which can be exploited in geometric quantization.

The local deformation theory for both  types of Lagrangian cycles
is well understood: its basic  fact is that every symplectic
manifold considered near some Lagrangian submanifold looks like
a small neighborhood of the zero section of the cotangent bundle
of the submanifold (the Darboux - Weinstein theorem).
Thus in this description the deformations of Bohr - Sommerfeld cycles
are presented by the graphs of exact forms while in the second case
the deformations are given by harmonic forms. This means that at least dimensionally
these two families of deformations are complement. 

On the space $\Cal L$ of all Lagrangian cycles one has the following bilinear forms.
First of all there is pure topological intersection form
$$
<S_1, S_2> = [S_1].[S_2] \in \Z
$$
where $[S_i] \in H_3(M, \Z)$. It's clear that it is skew symmetric.
At the same time for every two Lagrangian cycles intersecting
transversally, one can define the Floer homologies
$$
FH^*(S_1, S_2),
$$
([5]). For the pair $S_1, S_2$ we have a complex
$$
\de: \C^{S_1 \cap S_2} \to \C^{S_1 \cap S_2},
$$
defined by the set of the intersection points (we require transversality of the 
cycles). This construction results in a collection
of finite abelian groups labelled by integers and this collection
of the Floer homology groups gives us the following "bilinear"
form
$$
\theta(S_1, S_2) = \sum_{i=0}^3 (-1)^{i+1}rk FH^i(S_1, S_2),
$$
which is a symplectic analogue of $\chi$. 
If $S_1$ and $S_2$ are not transversal in $M$, then one can deform one of them
using some appropriate Hamiltonian deformation to establish the transversal 
picture and then apply the same arguments. For example, it's possible
to generalize the definition to the case
$$
S_1 = S_2 = S
$$
for a single cycle. Then in general the group concides with the usual de Rahm
cohomology of $S$ (at least dimensionally). All the details of
the long discussion can be found in [9] and other papers mentioned there.
For our story it's important that one can describe 
the deformation theory of the special Lagrangian cycles in a way similar
to the one described in LHS when we discussed the case of coherent sheaves:

1) the space $FH^1(S, S)$ is the space of infinitesimal deformations of
spLag cycle $S$;

2) the space $FH^2(S, S)$ is the space of obstructions;

3) the local model of the moduli space is given by some version of the  Kuranishi map

and

4) the spaces of the deformations
and the obstructions are dual and hence equidimensional.

As we've mentioned the deformation theory hardly depends on the topology of the cycle
itself and the topology of the embedding $S \to M$. But the real situation
is much more closer to the LHS: it was proved by McLean that all spLag deformations are
unobstructed. Together with the identification
$$
FH^1(S, S) = H^1(S, \R)
$$
(which takes place if the embedding $S \to M$ posesses some natural property which is
generic for the setup) it gives us that f.e. if $S$ is a homological sphere then
the moduli space of spLag cycles is presented by points. 

Now we are in position to formulate the Vafa conjecture which predicts more precise variant
of the mirror symmetry. For mirror partners $M, M'$ it should be a map
$$
mir: H^{2*}(M, \Z) \to H^3(M', \Z)
$$
induced by some one - to - one correspondence between $E$'s and $S$'s i.e. between
stable vector bundles over $M$ and spLag cycles inside the mirror partner
such that
$$
Ext^i(E_1, E_2) = FH^i(mir(E_1), mir(E_2)).
$$
One can draw this parallel further, considering the natural operations 
which can be applied for both objects: under the correspondence 
between $E$ and $S$ the operation of the Lagrangian connected sum
$$
S_1 \sharp S_2
$$
is dual to the extension operation:
$$
0 \to E_1 \to E \to E_2 \to 0
$$
for stable bundles $E_1$ and $E_2$.
 At the same time,  a notion
mentioned in LHS, is clarified: a stable vector bundle $E$ over
$M$ is called spherical (or by often used term "exeptional")
if
$$
Ext^1(E, E) = 0.
$$
By this strong version of the mirror conjecture,  $S = mir(E)$ has to be
a homological sphere 
$$
H^1(mir(E)) = H^2(mir(E)) = 0.
$$
Again it is reasonable to expect that the moduli space $\Cal M_S$ of spLag
realizations of this sphere is just a finite set of points:
$$
\Cal M_S = \{p_1, ..., p_d\}
$$
and the number of these points is a symplectic invariant
$$
SC_M(S) = \sharp(\sM_S).
$$
Then it would be natural to expect this number to be equal to
the Casson - Donaldson invariant of the stable bundle.

The symplectic part (or RHS) is not quite well understood 
yet. At the same time, the meaning of LHS is much more clear
and even suited for computations.  Here we will
present some method to deal with the LHS hoping that in the nearest future
matching investigations of RHS will appear.

\head 2. Degeneration method for the computation of the Casson - Donaldson
invariant
\endhead

\subhead 2.1. Degeneration principle
\endsubhead

One of the useful methods in algebraic geometry is based on
"degeneration principle": if one can reduce the situtation to
some appropriate degenerated case,  compute what is desired and
then prove that the number is invariant under the deformation,
then the problem is solved. As an example, we take an old standard problem
which was known 200 years ago: the question is how many lines
intersect 4 given skew ones in $\C \proj^3$. The regular construction
gives that there are exactly 2 such projective lines. But one can compute
this number using the degeneration principle as follows. Let us move
this couple of skew lines to the case when they are divided into two pairs
$(l_1, l_2), (m_1, m_2)$ such that 
$$
\aligned
<l_1, l_2> = \pi_l = \C\proj^2 \subset \C \proj^3 \\ 
<m_1, m_2> = \pi_m = \C \proj^2 \subset \C \proj^3\\
l_i \cap m_j = \emptyset. \\
\endaligned
$$
Then we see that the answer is 2: the first desired line is given
by the intersection points $l_1 \cap l_2$ and $m_1 \cap m_2$ while
the second one is given by the intersection $\pi_l \cap \pi_m$.

Now we would like to apply the degeneration principle to
the computation of the Casson - Donaldson invariant over Calabi
- Yau 3 - manifolds. Namely, let us solve the problem for some special
type of Calabi - Yau degeneration when a given 3 - manifold $M$ can be deformed
to a manifold
$$
M_0 = Y_+ \cup_S Y_-
\tag 2.1
$$
where $Y_{\pm}$ are smooth quasi Fano 3 - manifolds glued transversally
along a non singular surface $S$ of K3 - type which belongs simultaneously
to both the anticanonical systems:
$$
S \in \vert - K_{Y_{\pm}} \vert.
$$
Then if we compute the numbers in the framework of this
reduced picture and then prove that these numbers do not change under the
deformation, it will give us the answer in more general situation.

Let us set up the framework for the application of the degeneration principle
giving the following
\proclaim{Definition 2.1} Calabi - Yau 3 - manifold is called constructive
iff there exists an appropriate deformation of $M$ to $M_0$ which is
decomposed as in (2.1).
\endproclaim

The class of the constructive Calabi - Yau threefolds is wide enough:
of course, all complete intersections in every weighted projective space 
are constructive as well as any elliptic net, etc. Moreover,
one can observe that many rigid Calabi - Yau threefolds are constructive
(and it is a really amazing fact!). Here we place the following
\subheading{Example} The rigid Barth - Nieto - van Straten quintic
([10]) is the moduli space of abelian surfaces with the polarization of type
(2,6) and fixed theta - structure. But it could be realised using the following
model: consider the projective space $\C \proj^5$ equipped with homogeneous coordinates $(z_0, ..., z_5)$ and the corresponding system of the Newton
hypersurfaces
$$
S_k = \sum_{i=0}^5 z_i^k.
$$
Then the pencil of quintics 
$$
<S_5, S_2 \cdot S_3>
$$
in $\C \proj^4$ given by the linear equation $S_1 = 0$ contains
unique quintic with 130 nodes. It is the Barth - Nieto - van Straten
quintic.

To test the main idea of this article we will go in two different directions:
at the rest of Section 2 we will discuss how one can compute the invariant,
passing in the constructive case to some holomorphic symplectic setup
and reducing the computation to a usual computation in this setup.
In Section 3 we will show how one can construct "real" examples of the computation
when we will develop the corresponding gluing technique getting
constructive Calabi - Yau manifolds together with the results of
the computations.

So the next natural question  arises after the definition of the constructive Calabi - Yau manifolds is understood:
 what is the vector bundles over such reduced Calabi - Yau 3
- manifold $M_0$? A vector bundle $E$ on $M_0$ is a pair of vector bundles
$E_{\pm}$ over $Y_{\pm}$ such that their restrictions coincide:
$$
E_+|_S = E_-|_S.
$$
So the next step for us is to describe the geometry of vector bundles
over the flags of type $(S \subset Y)$ where $Y$ is a quasi Fano variety
(its definition see below) and $S$ is a K3 - surface from the anticanonical system.

\subhead 2.2. The geometry of the vector bundles on the flags
\endsubhead

We start with the following natural
\proclaim{Definition 2.2} A variety $Y$ is called a quasi Fano variety
if the anticanonical linear system contains a smooth K3 - surface and
$\chi(\Cal O_Y) = 1$.
\endproclaim

Well known examples of quasi Fano varieties are given by the blowing up
of the classical Fano varieties with centers on fixed surfaces of the anticanonical
systems.

Let $Y$ be a quasi Fano variety and $S$ be a fixed K3 - surface
from the anticanonical system. 
We will call such a pair  {\it a flag}.
This is a complex analogue of a real 3 - manifold with boundary.

For the quasi Fano varieties the arithmetical properties
of the Mukai lattice are slightly different from the CY - case, so
f.e. the bilinear form $\chi$ is not skew symmetric  and one has to decompose
it into the symmetric part and the skew - symmetric part denoting
these as $\chi_{\pm}$ respectively. Returning to the vector bundles 
on $Y$ we see that the symmetric form $\chi_+(E_1, E_2)$ depends
only on the first three components of the ring $H^{2*}(Y, \Z)$. 
Further, by the definition
$$
\Cal O_Y(K_X) = \Cal O_Y(-S).
$$
The canonical class (like any invertible sheaf) defines an automorphism
$T_{K_Y}$ of the Mukai lattice $L_Y$ equipped with the Mukai form
by the formula
$$
T_{K_Y}(m) = m \cdot e^{K_Y}.
$$
Restricting each vector bundle to the surface $S$ we obtain the following
map of the Mukai lattices:
$$
res: L_Y \to L_S = H^{2*}(S, \Z).
$$
The image of this map in the Mukai lattice of $S$ coincides with the image of
the following operator
$$
im(id - T_{-S}) \subset  L_M,
$$
where the last operator can be defined as follows
$$
\aligned
(id - T_{-S})(u_0, u_2, u_4, u_6) = (0, -K_Y \cdot u_0, u_2 \cdot S, u_4 \cdot S) = \\
(u_0, u_2 \cdot S, u_4 \cdot S) \in H^*(S).\\
\endaligned 
$$
Then the bilinear form
$$
<,> = \half res^* (,)
$$
is the preimage of the standard symmetric bilinear form on $H^{2*}(S, Z)$
given by
$$
(u, v) = - [v^* \cdot u]_4.
$$
Thus the restriction,  considered as a transformation, maps the root 
$$
\sqrt{td^+_Y}
$$
to the root of the Todd class of the K3 - surface
$$
(1, 0, 1) = \sqrt{td_S}
$$
since by the definition
of quasi Fano variety
$$
\frac{c_2(TY) \cdot K_Y}{24} = \chi(\Cal O_Y) = 1.
$$
The point is that for every Mukai vector $m$ corresponding to a vector bundle $E$
on $Y$  the Mukai vector of the restricted bundle is given by
$$
m(E|_S) = ch(E|_S) \cdot \sqrt{td_S} = (id - T_{-S})(m(E))
$$
over this K3 - surface. Moreover, the symmetric bilinear form
$\chi_+(E_1, E_2)$ on $Y$ is the lifting
of the classical symmetric form
$$
\chi(E_1|_S, E_2|_S)
$$
over $S$.
Consider  now "a geometric realization"
of $res$.  Namely,  since every vector bundle $E$ on $Y$ can be restricted
to $S$, we have a map
$$
r: \Cal M_E \to \Cal M_{E|_S}
$$
of the moduli spaces. The first result about these moduli spaces is purely arithmetical:
Proposition 11.2 of [13] ensures us that for any
simple vector bundle on $Y$ one has
$$
v. dim \Cal M_E = \half v. dim \Cal M_{E|_S}.
$$
We will see in a moment that this fact is much deeper
than just a purely numerical coincidence. 

The Bertini theorem gives us
$$
H^2(Y, \Z) = H^{1,1}(S) \cap H^2(S, \Z).
$$
A Mukai theorem says that for any primitive vector
$$
m = (u_0, u_2, u_4) \in L_S
$$
we have
$$
\{u_0>0; \quad u_2 \in H^{1,1}(S); \quad m^2 \geq -2\} \implies \sM_S(m) \neq \emptyset
$$
that is for this vector there exists a stable bundle $E_0$ over $S$ such that
for $m = m(E_0)$ the moduli space exists and
$$
dim \sM_S(m) \geq m^2 +2.
$$
Suppose that there exists a $-K_Y$ - stable bundle $E$ over $Y$ such that
$$
E_0 = E|_S.
$$
Then 
$$
\sM_Y(m) \neq \emptyset \quad {\text and} \quad dim \sM_Y(m) \geq \half m^2 + 1
$$
([8]).
To go further we need the following
\proclaim{Definition 2.3} 
1) A vector bundle $E$ on $Y$ is called regular if
$H^2(ad E) = 0$.

2) An irreducible component $\sM_Y(m)_0$ is called regular if a generic
bundle $E$ which belongs to this component is regular.
\endproclaim

These bundles are very important to us: for them we have
$$
dim \sM_Y(m)_0 = v.dim \sM_Y(m)_0 = \half m^2 + 1.
$$
Moreover, the restriction map
$$
r: \sM_Y(m)_0 \to \sM_S(m)
\tag 2.2
$$
is an immersion at any generic point. 
Both statements follow from the short exact cohomology sequence
$$
0 \to H^1(ad E) \to H^1(ad|_S) \to H^2(adE(K_Y)) \to 0
\tag 2.3
$$
constructed for the restriction sequence of $ad E$:
$$
0 \to ad E(K_Y) \to ad E \to ad E|_S \to 0
$$
(because $H^1(ad E(K_Y)$ is zero by the Serre duality). Moreover,
the continuation of the cohomological sequence gives
$$
0 \to H^2(ad E|S) \to H^3(ad E(K_Y)
$$
which ensures us that $H^2(adE|_S)$ is trivial if $E$ is simple:
$$
H^0(adE) = 0 \implies H^3(adE(K_Y)) = H^0(ad E)^* = 0.
$$
Therefore for simple bundles over $Y$ we have the following implication:
if $E$ is regular then  $E|_S$ is regular.
This gives us
\proclaim{Proposition 2.1} The restriction map $r$ (2.2)
is an embedding into the regular component $\sM_Y(m)_0$
such that
$$
dim \sM_S(m) = 2 dim \sM_Y(m)_0.
$$
\endproclaim

Now we invoke the  golomorphic symplectic geometry: 
every regular component $\sM_S(m)_0$ of vector bundles over
any K3 - surface admits the Mukai holomorphic symplectic structure
$$
\om_S: T\sM_S(m)_0 \to T^*\sM_S(m)_0.
\tag 2.4
$$
The crucial fact which underlies our further investigations is that
\proclaim{Proposition 2.2} The image of the restriction map
$$
r(\sM_Y(m)_0) \subset \sM_S(m)_0
$$
is a Lagrangian subvariety of $\sM_S(m)_0$ with respect to $\om_S$.
\endproclaim

The proof is very simple but quite illustrative.
The tangent space of $\sM_S(m)_0$ at any regular point $r(E)$ is isomorphic
to $H^1(ad|_S)$. Thus we can understand the second arrow in the exact sequence (2.3)
as the differential of the restriction map
$$
0 \to T_E \sM_Y(m)_0 @>dr>> T_{E|_S}\sM_S(m)_0
\tag 2.5
$$
and it is a monomorphism. On the other hand according to [8]
the Serre duality over $S$ is nothing but the restriction of the symplectic structure
$\om_S$ to the fiber of the tangent bundle that is
$$
T_{E|_S}\sM_S(m)_0 = H^1(ad E|_S) = H^1(ad E|_S)^* = T^*_{E|_S} \sM_S(m)_0
$$
and we can continue the sequence (2.5) dualizing it and using the identification (2.4)
which gives
$$
0 \to T_E \sM_Y(m)_0 @>dr>> T_{E|_S}\sM_S(m)_0 = T^*_{E|_S}\sM_S(m)_0
@>(dr)^*>> T^*_E \sM_Y(m)_0 \to 0.
\tag 2.6
$$
The sequence (2.6) is exact, as it is equivalent to (2.3), so 
$$
\om_S|_{T_E\sM_Y(m)_0} \cong 0
$$
and we are done.
Moreover, the same sequence (2.6) shows that in the regular case the normal
bundle should be isomorphic to the cotangent bundle, so
$$
N_{r(\sM_Y(m)_0) \subset \sM_S(m)_0} = T^* \sM_Y(m)_0.
\tag 2.7
$$

Proposition 2.2 suggests the introduction of a new integer invariant 
for Mukai vector $m$ over pair $(S, Y)$. Since $r(\sM_Y(m)_0)$ is 
a cycle of middle dimension in $\sM_S(m)_0$,  one can define its self intersection
index
$$
CD_{(S, Y)}(m) = [r(\sM_Y(m)_0)]^2
$$
which we call {\it the relative Casson - Donladson invariant}
of the pair $(S,  Y)$. In the compact and
non singular cases this number can be computed as the top Chern class of the
normal bundle:
$$
CD_{(S, Y)}(m) = c_{top}(N_{r(\sM_Y(m)_0), \sM_S(m)_0}),
$$
and (2.7) shows that in this case
$$
CD_{S, Y}(m) = \pm \chi (\sM_Y(m)_0).
\tag 2.8
$$
It is natural to exploit this relative version for the computations of
the absolute Casson - Donaldson invariant.

Thus the general strategy should be as follows:
for some Calabi - Yau threefold we are looking for
the deformation to a reducible "double" which is
glued from two quasi Fano varieties along the K3 - surface.
Then we reduce the question to the holomorphic symplectic
geometry over the K3 - surface $S$ where the moduli spaces
of vector bundles over the quasi Fanos live as 
holomorphic Lagrangian submanifolds.

At the same time this recipe can be adapted for computations
of some other type numbers: f.e. by the same procedure
we can compute the number of lines on generic $M_8$ which is 
a double cover of $\C \proj^3$ ramified at a  generic surface of degree 8
(we will study this case in details in the next section).
Then deforming this $M_8$ to a pair of $\C \proj^3_{\pm}$ glued along
a quartic $S$ we have two families of lines on each $\C \proj^3_{\pm}$
which are the Grassmannian $G(2,4)_{\pm}$. The intersection with
$S$ defines  maps
$$
r_{\pm}: G(2,4)_{\pm} \to Hilb^4(S).
$$
These maps obviously are embeddings. But the smooth variety $Hilb^4(S)$
has the Mukai holomorphic symplectic form $\om_S$. 
The images $r_{\pm}(G(2,4)_{\pm})$ are homotopy equivalent
smooth Lagrangian subvarieties. Hence the desired number is
$$
R_{M_8}(1) = 2 r(G(2,4))^2 = 2 c_{top}T^*G(2,4) = 12.
$$
Now we can ask: is it possible to compute some other
"classical" numbers using this degeneration method?
For example, let us compute the number $R_{Q_5}(k)$ of rational curves
of degree $k$ on a generic quintic $Q_5 \subset \C \proj^4$.
One can deform $Q_5$ to 
$$
M_0 = Q_2 \cup_S Q_3
$$
Then a rational curve $\ga$ of a degree, say, 10
degenerates to a reducible curve $\ga_3 \cup \ga_2$ where
$\ga_3$ is a cubic on $Q_3$ and $\ga_2$ is a conic on $Q_2$.
The intersection maps
$$
\aligned
r_+: \{\ga_2\} \to Hilb^6(S)\\
r_-: \{\ga_3 \} \to Hilb^6(S)\\
\endaligned
$$
are embeddings.
So what is the number
$$
\sharp[r_+(\{\ga_2\}) \cap r_-(\{\ga_3\})] = ?
$$
Note that the answer is unknown even for smaller
$k$s: the record is
$$
R_{Q_5}(5) = 229305888887625.
$$
However, the degeneration method gives some reason
why these numbers as cofficients of a generating function
are wrong (R. Pandharipande observed recently that the coefficient $n_{10}$
of the generating function doesn't give the number $R_{Q_5}(10)$).

In the next section we discuss the construction starting from the end:
this way we will find some particular examples of constructive Calabi - Yau
threefolds and study the vector bundles over these ones.

\head 3. Constructive Calabi - Yau threefolds
\endhead

\subhead 3.1. Deformations of  flags and  vector bundles
\endsubhead

The deformation theory of the pairs (K3 $\subset$ Fano) is quite similar to the
deformation theory for the complex manifold: one can construct over any flag 
a bundle $T(S, Y)$ (or a coherent sheaf) such that 

1) the space $H^1(T(S, Y))$ is the space of formal deformations;

2) the space $H^2(T(S, Y))$ is the space of obstructions;

3) and the corresponding Kuranishi map
$$
\Phi: H^1(T(S, Y)) \to H^2(T(S, Y))
$$
gives us a local model of the moduli space of deformations which is $\Phi^{-1}(0)$.
To construct the sheaf, consider the restriction sequence for the tangent bundle
$$
0 \to TY(-S) \to TY \to TY|_S
\tag 3.1
$$
together with the standard exact sequence on $S$
$$
0 \to TS \to TY|_S \to N_{(S, Y)} \to 0,
\tag 3.2
$$
where the last line bundle is the normal bundle to the surface in the threefold.
Let us combine two last epimorphisms from (3.1) and (3.2) getting
$$
TY \to N_{(S, Y)} \to 0
\tag 3.3
$$
and complete (3.3) to an exact sequence
$$
0 \to T(S, Y) \to TY \to N_{(S, Y)} \to 0.
$$
$T(S, Y)$ is the bundle (or the sheaf) which describes
the local deformation theory. Suppose that 
$$
H^1(S, N_{(S, Y)}) = 0.
$$
Then we can compare two Kuranishi maps, combining them in one diagram:
$$
\matrix 0 & \to & H^0(N_{S, Y}) & \to & H^1(T(S, Y)) & \to & H^1(TY) & \to & 0 \\
           &     & \downarrow   &     & \downarrow   &      & \downarrow &   &   \\
           &      & 0             & \to    & H^2(T(S, Y)) & \to    & H^2(TY) & \to & 0 \\
\endmatrix
$$

This gives us the following
\proclaim{Proposition 3.1} If the deformation of $Y$ is unobstructed then the same is true for the deformation of any pair $(S, Y)$.
\endproclaim

The definition of $T(S, Y)$ implies 
$$
0 \to TY(-S) \to T(S, Y) \to TS \to 0.
\tag 3.4
$$
Since our pair is of the type (K3  in Fano), the long cohomological sequence
gives for (3.4)
$$
\aligned
0 \to H^1(\La^2 \Om Y) \to H^1(T(S, Y)) \to H^1(\Om S)^* \\
\to H^1(\Om Y)^* \to H^2(T(S, Y)) \to 0 \\
\endaligned
$$
(here we use the equality
$$
E \otimes \La^3 E^* = \La^2 E^*
$$
for $rk E = 3$ - case together with the Serre duality). It's easy to see
that the homomorphism
$$
H^1(\Om S)^* \to H^1(\Om Y)^*
$$
from the last sequence is dual to the restriction map
$$
r: H^{1,1}(Y) \to H^{1,1}(S)
$$
(via the Dolbeaut isomorphism). This gives us
\proclaim{Proposition 3.2} The obstruction space is given by
$(ker\,\,\, r)^* \subset H^{2,2}(Y)$. In particular, if $Pic\,\, Y = \Z$, then
the deformation of $(S, Y)$ is unobstructed and 
according to Proposition 3.1 this implies that the deformation of $Y$ is unobstructed to.
\endproclaim

Indeed, the homomorphism
$$
H^1(\Om S)^* \to H^1(\Om Y)^*
$$
must be nontrivial. Thus it has to be an epimorphism.
Moreover, the space
$$
H^1(TY(K_Y)) = H^1(\La^2 \Om Y)
$$
is the space of the deformations of the pair which preserves the complex structure on $S$.

Recall that there exists a collection of
obstructions for the equivalence of n-th order thickening of our given K3 - surface $S$ in the quasi Fano variety $Y$ and its flat model. The first obstruction is given
by the class
$$
\om_1 \in H^1(TS \otimes N^*_{(S, Y)}).
$$
>From the standard exact sequence we have
$$
H^1(TS \otimes N^*_{(S, Y)}) = H^1(\La^2 \Om Y|_S).
$$
On the other hand, in our special case we have
$$
H^1(TS \otimes N^*_{(S, Y)}) = H^1(TS \otimes N_{(S, Y)})^*
$$
by the Serre duality. We will use these identifications for the "gluing procedure".

\subhead 3.2. Gluing procedure
\endsubhead

Starting with the configuration (2.1) we get over the fixed K3 - surface $S$ two normal bundles
$$
N_{S, Y_{\pm}} 
$$
which are completely different. The topological smoothing 
procedure is very similar to the topological surgery in the real case:
first of all, we cut a small neighborhood of the singular locus 
$S$ in $X_0$ considering small tubes
$$
S^1(N_{S, Y_{\pm}}) \subset N_{S, Y_{\pm}} 
$$
in the normal bundles which are
$$
N_{S, Y_{\pm}} = L_{-K_{Y_{\pm}}}.
$$
Removing small disc - bundles
$$
D^2(N_{S, Y_{\pm}})
$$
with the boundaries
$$
\partial D^2(N_{S, Y_{\pm}}) = - S^1 (N_{S, Y_{\pm}})
$$
from $Y_{\pm}$ we get some open singular threefold
$$
V_0 = D^2(N_{S, Y_+}) \cup_S D^2(N_{S, Y_-}).
$$
Thus one gets three real 6 - manifolds with boundaries:
$$
\aligned
V_0 \quad \quad | \quad \partial V_0 = S^1(N_{S, Y_+})
\cup S^1(N_{S, Y_-});\\
Y^0_{\pm} = Y_{\pm} - D^2(N_{S, Y_{\pm}}) \quad | \quad
\partial Y^0_{\pm} = S^1(N_{S, Y_{\pm}}); \\
\endaligned
$$
and the resulting $X_0$ is glued from these three pieces.

Now we can deform slightly the singular real 6 - manifold
$V_0$ preserving the boundary $\partial V_0$ by the following construction
([2]). Consider the following quadratic map of the bundles over $S$:
$$
q: N_{S, Y_+} \oplus N_{S, Y_-} \to N_{S, Y_+}
\otimes N_{S, Y_-} = L_{- K_{Y_+} - K_{Y_-}}.
$$
Any section $s \in H^0(L_{-K_{Y_+} - K_{Y_-}})$ can be regarded as an embedding
$$
i_s: S \to L_{-K_{Y_+} - K_{Y_-}}
$$
thus for any section we get a manifold
$$
V_s = q^{-1}(i_s(S)) \cap D^2(N_{S, Y_+}) \times_S D^2(N_{S,
Y_-}).
$$
 Choosing the neighborhoods small enough one can make the picture such that
the boundary becomes diffeomorphic  to
$$
\partial V_s = S^1(N_{S, Y_+}) \cup S^1(N_{S, Y_-}).
$$
Then we can glue this $V_s$ with $Y_{\pm}^0$ along the components of the boundaries and get 
a new compact real 6 - manifold $X_s$.
Moreover, if the zero set
$$
(s)_0 = C \in \vert -K_{Y_+} - K_{Y_-} \vert
$$
of the section $s$ is a smooth curve in $S$ then $V_s$ is non - singular
and  the construction gives us a topomodel of Calabi - Yau threefold.
If the curve $C$ admits some simple singularities then $V_s$ would be singular in these points but
applying the small resolution of these singular points one can get another topomodel of Calabi - Yau 
manifold of different topological type. It will be very usefull to get the complete list of the topomodels
which can be reached by this procedure.

Until now we discussed the gluing procedure from the point of view
of smooth real 6 - manifolds. But it is important that we can do this smoothing surgery preserving almost complex structures.

Now we go further describing the deformations of the complex structures over
the constructed topomodels. Recall ([2], [7]) that in our situation 
there exists a sheaf $T (Y_+, S, Y_-)$ which can be constructed 
in terms of $T(S, Y_{\pm})$ such that its first cohomology 
space presents the infinitesimal deformations of the reducible threefold
to reducible threefolds of the same topological type but the space 
$\Cal H^1$ of {\it all} infinitesimal deformations is more complicated:
it is included in the following exact sequence
$$
0 \to H^1(T(Y_+, S, Y_-)) \to \Cal H^1 \to H^0(N_{S, Y_+} \otimes N_{S,Y_-})
\to 0.
\tag 3.5
$$
However the obstruction space has precisely the same type
$$
\Cal H^2 = H^2(T(Y_+, S, Y_-))
$$
(5.1 of [2]). 
We can describe the sheaf $T(Y_+, S, Y_-)$ considering the entires
as disjoint flags $(S_{\pm} \subset Y_{\pm})$ together with the maps
$$
n_{\pm}: S_{\pm} \to S = Sing X_0.
$$
Constructing exact sequence (3.4) for every flag component
one gets the desired sheaf from the following exact sequence
$$
0 \to T(Y_+, S, Y_-) \to T(S_+, Y_+) \oplus T(S_-, Y_-) \to TS \to 0,
$$
where at the prefinal step we use $\half (n_+ + n_-)_*$.
Thus one has for $T(Y_+, S, Y_-)$  the exact sequence
$$
0 \to \oplus_{\pm} TY_{\pm}(-S) \to T(Y_+, S, Y_-) \to TS =
ker \half(n_+ + n_-)_* \to 0
$$
which induces the long cohomology sequence
$$
\aligned
0 \to \oplus_{\pm} H^1(\La^2 \Om Y_{\pm}) \to H^1(T(Y_+, S, Y_-))
\to H^1(\Om S)^* \\
\to \oplus_{\pm} H^1(\Om Y_{\pm})^*
\to H^2(T(Y_+, S, Y_-)) \to 0.\\
\endaligned
$$
The direct sum of the compositions $(r_{\pm} \cdot (n_{\pm})_*$ gives the map
$$
R^+_-=(r_+ \cdot(n_+)_*) \oplus (r_- \cdot(n_-)_*): H^{1,1}(Y_+)
\oplus H^{1,1}(Y_-) \to H^{1,1}(S).
$$
In terms of this map we can formulate the following
\proclaim{Proposition 3.3} The obstruction space is given by
$$
H^2(T(Y_+, S, Y_-)) = (ker R^+_-)^* \subset H^{1,1}(Y_+)^* \oplus
H^{1,1}(Y_-)^*.
$$
\endproclaim
Turning back to the definition of $\Cal H^1$ given by (3.5)
one can see that the deformation complex inducing the Kuranishi map
is 
$$
\aligned
0 \to H^1(T(Y_+, S, Y_-)) \to \Cal H^1 \to H^0(N_{S, Y_+} \otimes N_{S, Y_-})\\
@>\Psi>> H^2(T(Y_+, S, Y_-)) = (ker R^+_-)^* \subset H^{2,2}(Y_+)
\oplus H^{2,2}(Y_-).\\
\endaligned
$$
The precise description of $\Psi$ is contained in [2].

\proclaim{Corollary 3.1} 

1) Let $Pic Y_{\pm} = \Z$ and let $N_{S, Y_+} \otimes N_{S, Y_-}$
be generated by sections and nontrivial, then the dimension of the space $\Cal M_{X_0}$ of non singular deformations is
$$
v.dim \sM_{X_0} = h^{1,2}(Y_+) + h^{1,2}(Y_-) + h^0(N_{S, Y_+}
\otimes N_{S, Y_-}) - 1;
$$

2) if  $N_{S, Y_+} \otimes N_{S, Y_-} = \Cal O_S$ then the deformation
is unobstracted,
$$
v.dim \sM_{X_0} = h^{1,2}(Y_+) + h^{1,2}(Y_-) + 1
$$
and the body of the deformation family is smooth.
\endproclaim

The second statement is known from [2], [7]. The first statement can be proven as follows: consider the 1st order jet 
with the first obstruction
classes of the pairs $(S_{\pm}, Y_{\pm})$
$$
\om_1^{\pm} \in H^1(TS \otimes N^*_{S_{\pm}, Y_{\pm}}),
$$
described at the end of the previous subsection, and consider the natural
homomorphism
$$
\aligned
H^0(N_{S_+, Y_+} \otimes N_{S_-, Y_-}) \otimes H^1(TS \otimes N^*_{S_{\pm},
Y_{\pm}}) \\
@>c>> H^1(TS \otimes N_{S_{\pm}, Y_{\pm}}) = H^1(TS \otimes N^*_{S_{\pm}, Y_{\pm}})^*.\\
\endaligned
$$
It's easy to see that 1 - extension of the deformation given by a section
$s \in H^0(N_{S_+, Y_+} \otimes N_{S_-, Y_-})$ is constrained by only one
condition
$$
\om^-_1(c(s \otimes \om^+))=0.
$$
This implies the first statement.

Now we can illustrate the results as follows: suppose that some K3 - surface $S$
contained in a quasi Fano flag $S, Y$ admits an involution
$$
i: S \to S
$$
such that
$$
i^*(K_Y|_S) = K_Y|_S.
$$
Then take the reducible CY - threefold given by the gluing map with
$$
n_+ = id, \quad \quad n_- = i.
$$
and obtain
\proclaim{Corollary 3.2} Under the present conditions this double manifold
$X_0$ can be deformed to a smooth Calabi - Yau threefold.
\endproclaim

At the next subsection we add  the vector bundle ingredient
to our constructions and (since it's not too hard) consider
some particular examples.

\subhead 3.3. Vector bundles over constructive manifolds
\endsubhead.

The description of the vector bundles over the glued and deformed
threefold $X$ is very simple. If we denote the composition
$$
n_+ \cdot (n_-)^{-1}
$$
as $g$ ( thus $g$ is an automorphism of K3 - surface $S$) then
\proclaim{Proposition 3.4} A pair of stable vector bundles $E_{\pm}$
over the pair of quasi Fano threefolds $Y_{\pm}$ can be glued and deformed
to a vector bundle over $X$ if and only if
$$
E_+|_S = g^*(E_-|_S).
$$
\endproclaim

We may expect the infinitesimal rigidity for the vector bundles
on $X$. Our task is to describe the deformations of pairs
($X$ + a vector bundle). According to [2] this problem
can be reduced to the problem of the deformations for
the projectivizations of the vector bundle components.
But for this it's easy to consider only the deformations with smooth total
spaces. The criterium for this one is the following
\proclaim{Proposition 3.5} The total space of the deformations of
a reducible CY - threefold is smooth iff the tensor product of the normal bundles
$N_{S, Y_{\pm}}$ is trivial.
\endproclaim

Now we can reduce the general situation to the trivial tensor product case just
blowing up, say, the Fano variety $Y_+$ along the curve $C$. Thus 
we can expect that the construction can be performed  in a sufficiently
general case. 

 The simplest way to obtain a constructive threefold
is to take {\it double} of a flag. Let us add to
a flag $S \subset Y$ its copy and glue them along the fixed K3 - surface:
$$
2_S Y = (Y, S, Y)
$$
so
$$ 
n^{\pm} = id, \quad g = id.
$$
Certainly this double can be deformed to a smooth Calabi - Yau threefold
$X$ which type we will denote by the same symbol $2_SY$.
Any vector $m$ with non negative square from the Mukai lattice $M_Y$ defines the corresponding vector
$$
m \in M_{2_SY} \quad \quad {\text satisfies} \quad \quad m^2 = 0
$$
of the induced Mukai lattice of the double. For every regular component
$\sM_Y(m)_0$ of the stable vector bundle moduli space over $Y$ and for every
vector bundle $E \in \sM_Y(m)_0$ we can construct its {\it double} $2_SE$ 
 gluing two copies of $E$ along the restriction $E|_S$. As a result, instead of
the expected finite set of vector bundles we obtain some non transversal
component of the moduli space
$$
\sM_{2_SY}(m) = \{ 2_S E \quad | \quad E \in \sM_Y(m)_0 \}
$$
which obviosly has positive dimension being parametrized by $\sM_Y(m)_0$.
Thus from the first viewpoint the situation is not quite good.
But this large component decays to a finite set of vector bundles
when we deform the reducible double $2_S Y$ to a smooth $X$. The degree
of this finite set coincides with the number $CD_{S, Y}(m)$,
defined above in (2.8). The source of the coincidence is hidden in the following equality of two topological Euler characteristics
$$
\chi (\Cal D"(2_S E)) = \chi (\sM_Y(m)_0).
$$

At the rest of this section we place two examples of the construction
which are concentrated on objects well known in  algebraic geometry.

\subheading{Example 1} We consider for a smooth quartic surface $S$ in
$\C \proj^3$ the corresponding double of the flag. The resulting double
$2_S \C \proj^3$ is the degeneration of a smooth cover of
$\C \proj^3$ ramified along a smooth surface of degree 8.
For a vector bundle $E$ over $\C \proj^3$ of rank 2 with the Chern classes
$c_1 = 0, c_2 = 1$ its Mukai vector is presented by
$$
m(E) = 2 - \half P.D.(H);
$$
then the virtual dimension of the moduli space is given by
$$
v.dim \sM_{\C \proj^3} (2, 0, \half P.D.(H), 0) = m^2 + 1 = 5.
$$
It's well known that every stable bundle of the type is a mathematical instanton
given by a section of $\Om \C \proj^3(2)$ that is by a monad
$$
0 \to \Cal O_{\C\proj^3}(-1) @>s(-1)>> \Om \C \proj^3(-1) \to E \to 0.
$$
Therefore the compactification of the moduli space is constructed as
$$
\overline{\sM_{\C \proj^3}(2, 0, - \half P.D.(H), 0)} =
\C \proj^5 = \La^2 \C \proj^3.
$$
Now restricting every such vector bundle on our quartic
$S$ one gets another monad with the following display
$$
0 \to \Cal O_S(-1) @>s(-1)>> \Om \C \proj^3(-1)|_S \to E|_S \to 0.
$$
>From this it is easy to see that the restriction map
$$
res: \overline{\sM_{\C \proj^3}(2, 0, - \half P.D.(H),0)} \to 
\overline{\sM_S(2, 0, -2)}
$$
is an embedding. Thus we get
$$
CD_{(S, \C \proj^3}(2, 0, - \half P.D.(H), 0) = 6
$$
since it is the Euler characteristic of the projective space.
More geometrically, consider some general linear transformation
$$
g: \C \proj^3 \to \C \proj^3
$$
and the induced skew double $\C \proj^3 \cup_S g(\C \proj^3)$. Then a vector
bundle $E$ from the moduli space
$$
\overline{\sM_{\C \proj^3}(2, 0, - \half P.D.(H), 0)}
$$ 
defines the corresponding vector bundle on the skew double iff
$$
g^*(E) = E.
$$
Taking  $\La^2 \C \proj^3$ as a model of the compactification
consider the  transformation of the wedge square induced by $g$ and denote it
as $\La^2 g$. Then it's clear that
the condition $g^*(E) = E$ implies that $E$ corresponds to 
a fixed point of $\La^2 g$. Thus again we get that there are 6
vector bundles of this sort: the set $E_1, ..., E_6$ is just the six edges
of the simplex with the vertices in the fixed point of $g$ in $\C \proj^3$.

Now consider any smooth curve
$$
C = (s)_0, \quad s \in H^0(S, \Cal O_S(8))
$$
and blow it up as a curve in $\C \proj^3$, getting
$$
\si: \tilde{\C \proj^3} \to \C \proj^3.
$$
The gluing procedure then gives us an irreducible constructive Calabi - Yau threefold
$$
\tilde{\C \proj^3} \cup_S g (\C \proj^3)
$$
which can be deformed to a smooth threefold $X$ with smooth total deformation space. Then our six vector bundles give six
doubles
$$
\si^*(E_i) \cup_{S} E_i.
$$
One can prove that {\it the resulting vector bundles are infinitesimally rigid}.
It implies the fact that for any smooth double $X$ from the smooth deformation family
$$
CD_X((2, 0, \si^*(\half P.D(H)),0), (2, 0, -\half P.D.(H), 0)) = 6.
$$

\subheading{Example 2} We consider now the moduli space $MI_k$ of the mathematical
instantons, that is, the moduli space of the stable vector bundles over $\C \proj^3$
of rank 2 with the Chern classes $c_1 = 0, c_2 = k$ satisfying {\it the
instanton equation}
$$
h^1(E(-2)) = 0.
$$
Then, just like in the previous example, the restriction to $S$ must be an embedding by the same reason. 
The monad description shows that for general
linear transformation $g$ the induced action
$$
g^*: MI_k \to MI_k
$$
admits a finite set of fixed points
$$
E_1, ..., E_N
$$
which give the corresponding finite set of vector bundles over the double.
Note that any generic linear transformation defines a $\C^*$ -
action on $\C \proj^3$ and on the moduli space $MI_k$ as well. The computation
of the number $N$ (equal to the number of fixed points on $MI_k$) and of the
"instanton" component of the vector bundles over the double is
parallel to the computation of the Euler characteristic of $MI_k$ using
the Bott formula. Of course, this number $N$ is a term of CD - invariant of the double. However there
 exist some other components of the moduli space on $\C \proj^3$. Recall that any holomorphic vector bundle with topologically trivial    
determinant over any Fano variety of even index has the Atyiah - Rees invariant
$$
AR(E) = h^1(E(\half K_Y)) mod 2
$$
which distinguishes components of the moduli spaces. On the other hand, if 
$$
c_1(E) = 0
$$
then $E$ is anti self dual and the Serre duality induces the natural non degenerated skew symmetrical form on $H^1(E)$. Thus the equality
$$
h^1(E|_S) = 0 mod 2
$$
 always holds. Therefore for every mathematical instanton $E$ 
$$
AR(E) = 0
$$
and besides of $MI_k$ there exists another component $M_k$ (as a good example
we take $k = 3$). For a small $k$ it can be checked that the restriction
map embeds this component into the moduli space of
vector bundles over $S$. Thus
$$
res (MI_k) \cap res(M_k) = \emptyset.
$$
Hence  the Euler characteristics of the moduli spaces computed by the Bott
formula give the answer
$$
\aligned
CD_X((2, 0, - \si^*(\frac{5}{2} P.D.(H)), 0), (2, 0, -\frac{5}{2} P.D.(H), 0))
=\\
\chi(MI_3) + \chi(M_3).\\
\endaligned
$$

Some other examples are considered in [13].

\head Conclusion
\endhead

The author would propose a number of examples illustrating the algebraic geometry
of the constructions presented here. But let us emphasize just the questions
which are deeply important for the modern theoretical physics.

A standard question which is asked by physicists last time is about the possible topological types of 3 - dimensional Calabi - Yau manifolds. More precisely the question is just
on the number of these types. And even more concretely:
is this number bounded or not? The way proposed here could give an answer. Namely as it was already mentioned above if any 3 - dimensional Calabi - Yau
manifold is constructive then it were only finite number of
different topological types. Emphasize again that yet nobody knows
is this number finite or not.

On the other hand, one could ask an "adjoint" question:
{\it how many topological types  can be
derived by this construction of smooth deformations of threefolds}?
This question is quite natural  in the framework of algebraic geometry
(and much more simpler than the previos one).
But turning back one can ask a question of even more higher level:
{\it how many topological types of stable vector bundles can one get by this construction}?

At the same time the mirror conjecture dictates that our constructions
presented here should be compatible with some other aspects of mirror.
We mean that in the setup of Landau - Ginzburg models any Fano variety
admits a mirror partner. Very briefly, it is a pair
$$
((\C^*)^n, W_n),
$$
where $W_n$ is a function which is called {\it potential}.
If any 3 - dimensional Calabi - Yau manifold could be deformed to
a pair of flags then its mirror partner should be expressible in terms
of the Landau - Ginzburg models. Thus one can exploite the duality
to check the problem.

\enddocument